\documentclass[12pt,reqno]{amsart}
\usepackage{amsmath, amsthm, amsfonts, amscd,amsfonts,amssymb}

\newtheorem{thm}{Theorem}[section]
\newtheorem{cor}[thm]{Corollary}
\newtheorem{defn}[thm]{Definition}
\newtheorem{lem}[thm]{Lemma}

\usepackage[english]{babel}
\usepackage{ifthen}
\usepackage{graphics}
\usepackage{pifont}
\usepackage{float}
\usepackage{array}
\usepackage{hhline}
\usepackage{enumerate}
\usepackage{color}

\def\neg{{\rm neg}}

\usepackage{graphicx}
\usepackage{graphics}
\newcommand{\biindice}[3]%
{ #1\begin{array}[t]{c}
{\scriptstyle #2}\\
{\scriptstyle #3}\end{array}
}

\usepackage{pifont}
\usepackage{float}
\usepackage{array}
\usepackage{hhline}
\usepackage{enumerate}
\usepackage{color}

\author{Adel Hamdi}
\begin{document}
\title[crossings and nestings in permutations of type $B$]{Symmetric distribution of crossings and nestings in permutations of type $B$}
\maketitle
\begin{center}

\centerline{\footnotesize{Faculty of Science of Gabes, Department of Mathematics,}} \centerline{\footnotesize{ Cit\'e Erriadh 6072,  Zrig, Gabes, Tunisia}} \centerline{\footnotesize{and}}
\centerline{\footnotesize Universit\'e de Lyon, Universit\'e Lyon 1, Institut Camille Jordan,} \centerline{\footnotesize CNRS
UMR 5208, 43, boulevard du 11 novembre 1918}
\centerline{\footnotesize 69622, Villeurbanne Cedex, France}
\centerline{\footnotesize \texttt{aadel\_hamdi@yahoo.fr}}
\end{center}

\begin{abstract}
This note contains two results on the distribution of crossing numbers and nesting numbers in permutations of type B.
More precisely,
we prove a $B_n$-analogue of the symmetric distribution of crossings and nestings
of permutations due to Corteel (Adv. in Appl. Math., 38(2)(2007), 149-163) as well as the
symmetric distribution of $k$-crossings and $k$-nestings
of  permutations due to Burrill  et al. (DMTCS proc. AN, (2010), 461-468).

\end{abstract}
\section{ Introduction}
In the last years, many results on symmetric distributions of some statistics ``crossing" and ``nesting" have appeared in several combinatorial structures. At the heart of these results, on the set of matchings and partitions, there are Chen et al's and Kasraoui et al's
 theorems \cite{CDDS, KZ} on the interchanging crossing (resp., $2$-crossing, $k$-crossing) numbers and nesting (resp., $2$-nesting, $k$-nesting) numbers. Then, some extensions of type $B$ and $C$ have been given by Rubey and Stump~\cite{RS}, and Krattenthaler and de Mier~\cite{Kr,M} on the relation beteween increasing and decreasing chains in partitions and link partitions, and fillings of Ferrers shapes. On the set of permutations, Corteel~\cite{Co} has introduced the notion of crossings and nestings of permutations and proved that for any fixed number of weak exceedances, the distribution of crossing numbers and nestings numbers of permutations is symmetric. Recently, Burrill et al \cite{BMP} have proved a similar result for $k$-crossings and $k$-nestings of permutations. The purpose of this paper is to extend the last two results to their analogue of type $B$.
\section{ Definitions and main results}
For a positive integer $n$, let [$n$] := $\{1, 2, \ldots, n\}$. A type $B$ permutation of rank $n$ is an integer sequence
 $\sigma:= (\sigma(1), \sigma(2), \ldots,\sigma(n))$ such that $\{|\sigma(1)|, ...,|\sigma(n)|\} = [n]$.
  In this paper, we shall identify $\sigma$ with a permutation of [$-n$, $n$] := $\{-n, \ldots, -2, -1,$ $1, \ldots, n\}$ by $\sigma$($-i$) = $-\sigma(i)$ for each $i$ $\in$ [$n$]. Let $\neg(\sigma)$ be the \textit{number of negative numbers} in $\{\sigma(1), \dots,\sigma(n)\}$, and $B_n$ the set of type $B$  permutations of rank $n$.

\hspace{0.2cm} In the sequel, we use the natural order of integers in  $\mathbb{Z}$.

\hspace{0.2cm}As in \cite{Co}, it is convenient to represent a permutation $\sigma$ $\in$ $B_n$ by a \emph{permutation diagram} $G=(V,E)$, where $V=[-n,\;n]$ is the
vertex set, and $E$ is the set of edges $(i,\sigma(i))$ for $i\in [-n,\;n]$ such that
the vertices $-n$, \ldots, $-2$, $-1$, 1, 2, \ldots, $n$ are  arranged from left to right on a straight line. We draw an arc from $i$ to $\sigma(i)$ above (resp. under) the line if $i \leqslant \sigma(i)$ (resp.  otherwise) such that two arcs cross at most once.
A  permutation diagram is given in Fig.~1.

\begin{picture}(4,2.8)(-4,-1.3)
\qbezier(.4,.4)(1.2,1.2)(2,.4)
\qbezier(2.4,.4)(3.,1.2)(3.6,.4)
\qbezier(3.6,.4)(3.8,.8)(4,.4)
\qbezier(.4,.4)(.6,-.04)(.8,.4)
\qbezier(.8,.4)(1.4,-.55)(2,.4)
\qbezier(2.4,.4)(3.2,-.6)(4,.4)
\qbezier(0,.4)(1.5,1.5)(2.8,.4)
\qbezier(1.6,.4)(3.1,-.9)(4.4,.4)
\qbezier(1.6,.4)(3.1,1.5)(4.4,.4)
\qbezier(0,.4)(1.5,-.9)(2.8,.4)
\qbezier(1.2,.4)(1.4,.58)(1.2,.65)
\qbezier(1.2,.4)(1.,.58)(1.2,.65)

\qbezier(3.2,.4)(3.4,.58)(3.2,.65)
\qbezier(3.2,.4)(3.,.58)(3.2,.65)

\put(1.25,.32){\scriptsize{-3}}
 \put(.85,.32){\scriptsize{-4}}
\put(1.65,.32){\scriptsize{-2}}
 \put(2.05,.32){\scriptsize{-1}}
 \put(2.45,.32){\scriptsize{1}}
 \put(2.82,.32){\scriptsize{2}}
 \put(3.24,.32){\scriptsize{3}}
 \put(3.65,.32){\scriptsize{4}}
\put(4.04,.32){\scriptsize{5}}
\put(.46,.32){\scriptsize{-5}}
\put(4.45,.32){\scriptsize{6}}
\put(.06,.32){\scriptsize{-6}}

 \put(2.4,.4){\circle*{0.1}}
 \put(.4,.4){\circle*{0.1}}
\put(1.6,.4){\circle*{0.1}}
 \put(2,.4){\circle*{0.1}}
 \put(1.2,.4){\circle*{0.1}}
 \put(2.8,.4){\circle*{0.1}}
 \put(3.2,.4){\circle*{0.1}}
 \put(3.6,.4){\circle*{0.1}}
\put(4,.4){\circle*{0.1}}
\put(.8,.4){\circle*{0.1}}
 \put(4.4,.4){\circle*{0.1}}
 \put(0,.4){\circle*{0.1}}
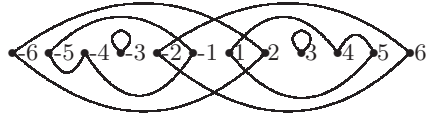
\put(-1,-.8){\footnotesize{Fig. 1. The permutation diagram of $\sigma= (4, -6, 3, 5, 1, -2)$}.}
\end{picture}

\hspace{0.2cm}We call the set of arcs that are above (resp. under) the line the upper (resp. under) permutation diagram and denoted
\textit{Upp}($\sigma$) (resp. \textit{Und}($\sigma$)).

\hspace{0.2cm}We start with an easy lemma that follows immediately from the definition of the permutation diagram since there is an easy bijective between upper and under diagarms.

\begin{lem}
Let $\sigma$ $\in$ $B_n$. The diagram of $\sigma$ is completely determined by the Upp($\sigma$).
\end{lem}

\hspace{0.2cm}Note that there are five geometric patterns for two arcs above the line as illustrated in Fig.~2 by

\begin{picture}(4,1.5)(-5,-.4)
\qbezier(-3,.4)(-3.4,1.)(-3.8,.4)
\qbezier(-3.4,.4)(-3.8,1.)(-4.2,.4)
\put(-3.4,.4){\circle*{0.1}}
\put(-3.8,.4){\circle*{0.1}}
\put(-4.2,.4){\circle*{0.1}}
\put(-3,.4){\circle*{0.1}}
\put(-3.7,-.15){(i)}

\qbezier(-.8,.4)(-1.,.8)(-1.2,.4)
\qbezier(-1.2,.4)(-1.4,.8)(-1.6,.4)
\put(-1.2,.4){\circle*{0.1}}
\put(-1.6,.4){\circle*{0.1}}
\put(-.8,.4){\circle*{0.1}}
\put(-1.5,-.15){(ii)}

\qbezier(.6,.4)(1.2,1.4)(1.8,.4)
\qbezier(1.,.4)(1.2,.8)(1.4,.4)
\put(1.,.4){\circle*{0.1}}
\put(1.4,.4){\circle*{0.1}}
\put(1.8,.4){\circle*{0.1}}
\put(.6,.4){\circle*{0.1}}
\put(.9,-.15){(iii)}

\qbezier(3.2,.4)(3.8,1.4)(4.4,.4)
\qbezier(3.8,.4)(4,.58)(3.8,.65)
\qbezier(3.8,.4)(3.6,.58)(3.8,.65)
\put(3.2,.4){\circle*{0.1}}
\put(3.8,.4){\circle*{0.1}}
\put(4.4,.4){\circle*{0.1}}
\put(3.6,-.15){(iv)}

\qbezier(5.8,.4)(6,.8)(6.2,.4)
\qbezier(6.6,.4)(6.8,.8)(7,.4)
\put(6.2,.4){\circle*{0.1}}
\put(6.6,.4){\circle*{0.1}}
\put(7,.4){\circle*{0.1}}
\put(5.8,.4){\circle*{0.1}}
\put(6.1,-.15){(v)}
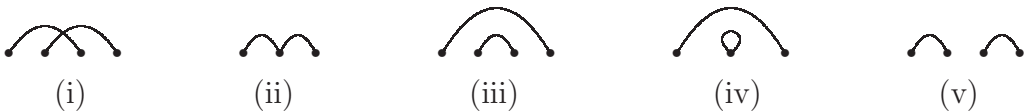
\put(-2.7,-.8){\footnotesize{Fig. 2. Five patterns between two arcs above the line.}}
\end{picture}\\\\

\hspace{0.2cm}These patterns are called: (i) a proper crossing, (ii) a skew crossing, (iii) a proper nesting, (iv) a skew nesting and (v) an alignment.
In another sense, one can recover these geometric patterns as in the two following definitions.

\hspace{0.2cm}The first is the notion of crossings of type $B$ given by Corteel et al. in \cite{CMW} as follows.

\begin{defn}
Let $\sigma\in$ $B_n$. The number of weak exceedances of $\sigma$, denoted by
$wex_B(\sigma)$, is the cardinality of the set $\{j \in[n]; \sigma(j) \geq j\}$. For two integers $i$ and $j$ in $[n]$,
 two arcs $(i, \sigma(i))$ and $(j, \sigma(j))$  form a crossing of $\sigma$ if they satisfy either the relation $  i < j \leqslant \sigma(i) < \sigma(j)$
(upper crosing), or $  -i < j \leqslant \sigma(-i) < \sigma(j)$ (upper crosing) or $  \sigma(i) < \sigma(j) < i < j $ (lower crosing).
\end{defn}
\hspace{0.2cm}Similarly, in the second, we can define the notion of nesting of type $B$.
\begin{defn}
Let $\sigma$ $\in$ $B_n$. A pair of arcs $($i, $\sigma(i))$ and $($j, $\sigma(j))$, with $i$ and $j$ in $[$n$]$,  is a nesting of $\sigma$ if they satisfy either
the relation $ i < j \leqslant \sigma(j) < \sigma(i)$ (uppeder nesting), or $ -i < j \leqslant \sigma(-j) < \sigma(i)$ (upper nesting) or $ \sigma(j) < \sigma(i) < i < j  $ (lower nesting).
The number of crossings $($resp., nestings$)$ of $\sigma$ is denoted by $cro_B(\sigma)$ $($resp., $nes_B(\sigma))$.
\end{defn}
\textbf{Example 1.} Let $\sigma$ = (4, $-6$, 3, 5, 1, $-2$) $\in$ $B_6$. Then the nestings in $\sigma$ are $\{(-2$, $\sigma(-2)$),
(1, $\sigma(1)$)$\}$, $\{$(3, $\sigma(3)$), (1, $\sigma(1)$)$\}$, $\{$(3, $\sigma(3)$), ($-2$, $\sigma(-2)$)$\}$, $\{$(4, $\sigma(4)$), ($-2$, $\sigma(-2)$)$\}$ and $\{$(6, $\sigma(6)$), (5, $\sigma(5)$)$\}$.
The crossings are $\{(-6$, $\sigma(-6)$), (1,$\sigma(1)$)$\}$, $\{$(1, $\sigma(1)$), (4,$\sigma(4)$)$\}$, $\{$(5, $\sigma(5)$), (2,$\sigma(2)$)$\}$ and $\{$(6, $\sigma(6)$), (2,$\sigma(2)$)$\}$ (see Fig.~1).
Hence $nes_B(\sigma)$ = 5 and $cro_B(\sigma)$ = 4.

\hspace{0.2cm}The following is our $B_n$-analogue of Corteel's result for type A permutations \cite[Proposition 4]{Co}.
\begin{thm}
The number of permutations in $B_n$ with $k$ weak exceedances, $l$ minus signs, $i$ crossings and $j$
nestings is equal to the number of permutations in $B_n$ with $k$ weak exceedances, $l$ minus signs, $i$ nestings and $j$ crossings.
In other words, we have
\begin{align}
\sum_{\sigma\in B_n}p^{nes_B(\sigma)}q^{cro_B(\sigma)}y^{wex_B(\sigma)}a^{neg(\sigma)}
= \sum_{\sigma\in B_n}p^{cro_B(\sigma)}q^{nes_B(\sigma)}y^{wex_B(\sigma)}a^{neg(\sigma)}.
 \end{align}
 \end{thm}
\hspace{0.2cm}Note that the $a$ = 0 case of Theorem 2.4 corresponds to Proposition 4 in \cite{Co}.

\hspace{0.2cm}Now, we extend the definition of $k$-crossings and $k$-nestings for permutations of type $A$ in \cite{BMP} to permutations of type $B$.
\begin{defn}
Let $\sigma$ $\in$ $B_n$. A set $\{ a_1, a_2, \ldots, a_k\}$ of k integers in $[n]$ is a $k$-crossing of $\sigma$ if they satisfy either the relation
 $a_1 < a_2 <\ldots< a_k \leq \sigma(a_1) < \sigma(a_2) < \ldots < \sigma(a_k)$ $($upper $k$-crossing$)$,
or $-a_1 < a_2 <\ldots< a_k \leq $ $-\sigma(a_1) < \sigma(a_2) < \ldots < \sigma(a_k)$ $($upper $k$-crossing$)$ or
$\sigma(a_k) < \sigma(a_{k-1}) < \ldots < \sigma(a_1) < a_k < a_{k-1} <\ldots< a_1$ $($lower $k$-crossing$)$.
\end{defn}

\begin{defn}
Let $\sigma$ $\in$ $B_n$. A set $\{ a_1, a_2, \ldots, a_k\}$ of k integers in $[n]$ is a $k$-nesting of $\sigma$ if they satisfy either the relation
$a_1 < a_2 <...< a_k \leq \sigma(a_k) < \sigma(a_{k-1}) < \ldots < \sigma(a_1)$ $($upper $k$-nesting$)$,
or $-a_1 < a_2 <\ldots< a_k \leq$ $-\sigma(a_k) < \sigma(a_{k-1}) < \ldots< \sigma(a_1)$ $($upper $k$-nesting$)$ or
 $\sigma(a_k) < \sigma(a_{k-1}) < \ldots < \sigma(a_1) < a_1 < a_{2} <\ldots< a_k$ (lower $k$-nesting).
\end{defn}

\hspace{0.2cm}As in \cite{BMP}, the \textit{k-crossing number} (resp., \textit{k-nesting number}) of a permutation $\sigma$ of type $B$, denoted by $cro^*_B(\sigma)$ (resp., $nes^*_B(\sigma)$) is the size of the largest $k$ such that $\sigma$ contains a $k$-crossing (rsep. $k$-nesting).\\\\
\textbf{Example 2.} Let $\sigma$ = (4, 5, 6, 2, $-3$, $-1$) $\in$ $B_6$. Then we have $cro^*_B(\sigma)$ = 4 and $nes^*_B(\sigma)$ = 2 that are illustrated respectively, in Fig.~3,
by $\{5, 1, 2, 3\}$ and $\{$4, 5$\}$ or $\{$4, 6$\}$ since $-5 < 1 <2 < 3 \leqslant -\sigma(5) < \sigma(1) < \sigma(2) < \sigma(3)$, $\sigma(5) < \sigma(4) < 4 < 5$
and $\sigma(6) < \sigma(4) < 4 < 6$.

\begin{picture}(-4,2.56)(-4,-1.2)
\qbezier(.4,.4)(1.9,1.75)(3.2,.4)
\qbezier(.8,.4)(1.2,1.)(1.6,.4)
\qbezier(1.2,.4)(2.6,-1.05)(4,.4)

\qbezier(2,.4)(3.2,-.9)(4.4,.4)

\qbezier(2.4,.4)(3,1.3)(3.6,.4)
\qbezier(2.8,.4)(3.4,1.3)(4,.4)
\qbezier(3.2,.4)(3.8,1.3)(4.4,.4)

\qbezier(2.8,.4)(3.2,-.2)(3.6,.4)

\qbezier(.,.4)(.6,-.4)(1.2,.4)
\qbezier(.4,.4)(1,-.4)(1.6,.4)
\qbezier(.8,.4)(1.4,-.4)(2,.4)

\qbezier(0,.4)(1.,1.6)(2.4,.4)

\put(1.25,.32){\scriptsize{-3}}
 \put(.85,.32){\scriptsize{-4}}
\put(1.64,.32){\scriptsize{-2}}
 \put(2.05,.32){\scriptsize{-1}}
 \put(2.45,.32){\scriptsize{1}}
 \put(2.85,.32){\scriptsize{2}}
 \put(3.27,.32){\scriptsize{3}}
 \put(3.65,.32){\scriptsize{4}}
\put(4.04,.32){\scriptsize{5}}
\put(.45,.32){\scriptsize{-5}}
\put(4.45,.32){\scriptsize{6}}
\put(.04,.32){\scriptsize{-6}}

 \put(2.4,.4){\circle*{0.1}}
 \put(.4,.4){\circle*{0.1}}
\put(1.6,.4){\circle*{0.1}}
 \put(2,.4){\circle*{0.1}}
 \put(1.2,.4){\circle*{0.1}}
 \put(2.8,.4){\circle*{0.1}}
 \put(3.2,.4){\circle*{0.1}}
 \put(3.6,.4){\circle*{0.1}}
\put(4,.4){\circle*{0.1}}
\put(.8,.4){\circle*{0.1}}
 \put(4.4,.4){\circle*{0.1}}
 \put(0,.4){\circle*{0.1}}
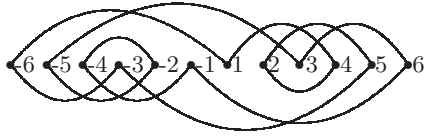
 \put(-1.2,-1){\footnotesize{Fig. 3. The permutation diagram of $\sigma$ = (4, 5, 6, 2, -3, -1)}.}
\end{picture}

\hspace{0.2cm}We also recall the following definition from \cite{BMP}.
Let $\sigma$ $\in$ $B_n$. The \emph{degree sequence} of the upper permutation diagram $Upp(\sigma)$ is the sequence $(indegree_{\sigma}(i), outdegree_{\sigma}(i))_{i \in[-n,n]}$,
where $indegree_{\sigma}(i)$ (resp. $outdegree_{\sigma}(i)$) is the left (resp. right) degree of the vertex $i$, $i.e.$, the number of arcs joining $i$ to a vertex $j$ with
$j<i$ (resp. $j>i$). If an upper permutation diagram $Upp(\sigma)$ has $d$ as its degree sequence (some other sources call this left-right degree sequence),
we say that $Upp(\sigma)$ is a diagram on $d$.\\
But there is a straightforward difference that we do not put a loop (an arc if $\sigma(i) = i$) on the isolated vertex $i$ with negative index in $Upp(\sigma)$,
 $i.e.$, we put (0,0) as a degree. By Lemma 2.1, we limit ourselves to study the upper permutation diagram of type B. The vertices with degree (0,1) (resp., (1,0), (1,1)) are called $openers$
(resp., $closers$, closer-opener or transient). For instance, if we let $\sigma$ = (4, 5, 6, 2, -3, -1), then the degree sequence of the upper permutation diagram of $\sigma$ is
$$
d:=d(\sigma) = (0,1)(0,1)(0,1)(0,0)(1,0)(0,0)(1,1)(0,1)(1,1)(1,0)(1,0)(1,0).
$$
Let $B_n^d$ be the set of the permutations in $B_n$ that has the degree sequence d.

\hspace{0.2cm}The following is our $B_n$-analogue of \cite[Theorem 1]{BMP}.
\begin{thm}
Let $NC_{B_n^d}(i,j,m)$ be the number of permutations in $B_n$ with i-crossings, j-nestings, $m$ minus signs and degree sequence specified by $d$. Then
\begin{align}
NC_{B_n^d}(i,j,m) = NC_{B_n^d}(j,i,m).
\end{align} \textit{In other words, we have} \\
\begin{align}
\sum_{\sigma\in B_n^d}x^{nes^*_B(\sigma)}y^{cro^*_B(\sigma)}z^{neg(\sigma)} = \sum_{\sigma\in B_n^d}p^{cro^*_B(\sigma)}q^{nes^*_B(\sigma)}z^{neg(\sigma)}.
 \end{align}
\end{thm}
\hspace{0.2cm}Note that when $z$ = 0 in (3), we recover Theorem 1 of \cite{BMP}.

\hspace{0.2cm}Now, we sketch the opener, closer and transient vertices of the upper (resp., under) permutation diagram with degree (0,1), (1,0) and (1,1)
 (resp., (1,0), (0,1), (1,1)) respectively. Then a vertex is said to be:

\begin{picture}(-4,1.9)(-5,-.45)
\qbezier(1.4,.9)(1.55,1.15)(1.7,1.2)
\put(1.4,.35){\circle*{0.1}}
\qbezier(3.35,.9)(3.2,0.62)(3.05,.6)
\put(3.35,.9){\circle*{0.1}}
\put(2,.85){(resp.,\;\;\;\;\;\;\;\;\;\ )}
\put(-4.5,.85){(i)\;\; an opener if it is illustrated by}

\qbezier(1.4,.4)(1.25,.65)(1.1,.7)
\put(1.4,.85){\circle*{0.1}}
\qbezier(3.35,.4)(3.43,0.25)(3.7,.18)
\put(3.35,.45){\circle*{0.1}}
\put(2,.35){(resp.,\;\;\;\;\;\;\;\;\;\ )}
\put(-4.5,.35){(ii)\; a closer if it is illustrated by}

\qbezier(1.4,.)(1.55,.25)(1.7,.3)
\qbezier(1.4,.)(1.25,.25)(1.1,.3)
\put(1.4,-.05){\circle*{0.1}}
\qbezier(3.35,.)(3.55,-0.3)(3.7,-.33)
\qbezier(3.35,.)(3.15,-0.3)(3.05,-.33)
\put(3.35,-.0){\circle*{0.1}}
\put(2,-.05){(resp.,\;\;\;\;\;\;\;\;\;\ ).}
\put(-4.5,-.05){(iii) a transient if it is illustrated by}
\end{picture}

\hspace{0.2cm}We shall prove Theorem 2.4  in Section 3 by constructing an explicit involution on $B_n$ that interchanges the number of crossings and number of nestings.
In fact, it is an extension of the involution defined in \cite{KZ}. To prove the Theorem 2.7  in  Section~4, we shall adopt the map
defined by de Mier in \cite{M} to $B_n$.
\section{Proof of Theorem 2.4}
First, for each $\sigma$ $\in$ $B_n$, the number of crossings of $\sigma$ is equal, in $Upp(\sigma)$, to the number of proper crossings plus the number of
transient vertices $i$ with $i$ in [$n$]. Similarly, the number of nestings of $\sigma$ is equal, in $Upp(\sigma)$, to the number of proper nestings plus
 the number of arcs ($i$, $k$) with $i$ and $k$ in [$-n, n$] such that there exist two fixed vertices $j$ and $-j$ satisfy $i<|j|<k$. For instance, in the left diagram of Fig.~4, we can count the 4 crossings and the 5 nestings of $\sigma$ given in Fig. 1.

\hspace{0.2cm}Now, we introduce some notations. For a positive integer $n$, let $\Lambda_n$ be the set of the subsets of [$n$] and $\overline{B}_n$ the set of $Upp$($\sigma$) for each $\sigma$ $\in$ $B_n$. Let also $F$ and $T$ be two maps defined by: for each $\sigma$ $\in$ $B_n$, $F(\sigma) := \{ j \in [n]; \sigma(j) = j\}$ and  $T(\sigma):= \{ j \in [n]; \;\ j\;\ is\,\ an\;\ upper\;\ transient\;\ vertex\;\ of\;\ \sigma \}$. We notice that, in \cite{RS}, Rubey and Stump have studied the symmetry distribution of the number of crossings
 and number of nestings in a kind of set partitions of type B. Then, we study a simile result in $B_n$ where we count the transient vertex
(resp. an arc covers a fixed vertex) as a crossing (resp. nesting) and our arrangement of the set [$-n, n$] is different of their.

\textit{Proof of Theorem 2.4}.
There are two steps.

\textit{ First step.}
 Let $\sigma$ $\in$ $B_n$. We define a map $\psi$ by: let ($Upp(\sigma$), $F(\sigma$), $T(\sigma$)) $\in$ $\overline{B}_n \times \Lambda^2_n$.
Then $\psi$ transforms each element $i$ of $F(\sigma)$ to an arc ($i$, $i'$) (see (a)) and each element $j$ of $T(\sigma)$ to a proper crossing (see (b)),
 which in the two cases, we have $i<i'$ (resp. $j<j'$) and no vertex between $i$ and $i'$ (resp. $j$ and $j'$). We adopte from \cite{BMP} the following graphs\\
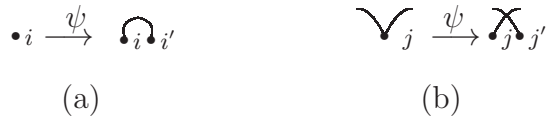
\begin{picture}(1,1.2)(-3,.2)
\qbezier(4.4,.9)(4.55,1.2)(4.7,1.2)\put(4.6,.8){\scriptsize{$j$}}
\qbezier(4.4,.9)(4.25,1.2)(4.1,1.2)
\put(4.4,.9){\circle*{0.1}}
\put(4.9,.8){$\longrightarrow$}\put(5.05,1.){$\psi$}
\qbezier(5.6,.9)(5.75,1.2)(5.9,1.2)
\qbezier(5.9,.9)(5.75,1.2)(5.6,1.2)
\put(5.7,.8){\scriptsize{$j$}}\put(6.,.8){\scriptsize{$j'$}}
\put(5.6,.9){\circle*{0.1}}\put(5.9,.9){\circle*{0.1}}
\put(4.8,.1){(b)}
\end{picture}
\begin{picture}(1,1.2)(5,.6)
\put(7.2,1.3){\circle*{0.1}}
\put(7.5,1.2){$\longrightarrow$}\put(7.75,1.4){$\psi$}
\qbezier(8.4,1.3)(8.35,1.46)(8.55,1.5)\put(7.3,1.2){\scriptsize{$i$}}
\qbezier(8.7,1.3)(8.75,1.46)(8.55,1.5)
\put(8.5,1.15){\scriptsize{$i$}}\put(8.8,1.15){\scriptsize{$i'$}}
\put(8.7,1.26){\circle*{0.1}}\put(8.4,1.26){\circle*{0.1}}
\put(7.7,.5){(a)}
\end{picture}

 For instance, the permutation diagram of $\sigma$ in Example 1 has a fixed vertex indexed by 3 and a transient vertex indexed by 4,
see the left diagram of Fig.~4. Then, its inverse $\psi^{-1}$ reduces each two vertices $i$ and $i'$ introduced by $\psi$ into a one vertex $i$.

\begin{picture}(4,1.5)(-1,-.7)
\qbezier(.4,.4)(1.2,1.2)(2,.4)
\qbezier(2.4,.4)(3.,1.2)(3.6,.4)
\qbezier(3.6,.4)(3.8,.8)(4,.4)

\qbezier(0,.4)(1.5,1.5)(2.8,.4)
\qbezier(1.6,.4)(3.1,1.5)(4.4,.4)

\put(1.25,.32){\scriptsize{-3}}
\put(.85,.32){\scriptsize{-4}}
\put(1.65,.32){\scriptsize{-2}}
\put(2.05,.32){\scriptsize{-1}}
 \put(2.47,.32){\scriptsize{1}}
 \put(2.85,.32){\scriptsize{2}}
 \put(3.24,.32){\scriptsize{3}}
 \put(3.68,.32){\scriptsize{4}}
\put(4.06,.32){\scriptsize{5}}
\put(.45,.32){\scriptsize{-5}}
\put(4.47,.32){\scriptsize{6}}
\put(.04,.32){\scriptsize{-6}}

 \put(2.4,.4){\circle*{0.1}}
 \put(.4,.4){\circle*{0.1}}
\put(1.6,.4){\circle*{0.1}}
 \put(2,.4){\circle*{0.1}}
 \put(1.2,.4){\circle*{0.1}}
 \put(2.8,.4){\circle*{0.1}}
 \put(3.2,.4){\circle*{0.1}}
 \put(3.6,.4){\circle*{0.1}}
\put(4,.4){\circle*{0.1}}
\put(.8,.4){\circle*{0.1}}
 \put(4.4,.4){\circle*{0.1}}
 \put(0,.4){\circle*{0.1}}

\end{picture}
\begin{picture}(4,2.)(-3.,-.7)

\qbezier(.4,.4)(1.2,1.2)(2,.4)
\qbezier(2.4,.4)(3.035,1.2)(3.67,.4)
\qbezier(3.51,.4)(3.755,.8)(4,.4)
\qbezier(0,.4)(1.4,1.5)(2.8,.4)
\qbezier(1.6,.4)(3.,1.5)(4.4,.4)

\qbezier(3.13,.4)(3.04,.55)(3.2,.58)
\qbezier(3.27,.4)(3.33,.55)(3.2,.58)

\put(1.25,.32){\scriptsize{-3}}
 \put(.85,.32){\scriptsize{-4}}
\put(1.65,.32){\scriptsize{-2}}
 \put(2.05,.32){\scriptsize{-1}}
 \put(2.475,.32){\scriptsize{1}}
 \put(2.85,.32){\scriptsize{2}}
 \put(3.06,.05){\scriptsize{3}}\put(3.19,.05){\scriptsize{3'}}
 \put(3.45,.05){\scriptsize{4}}\put(3.6,.05){\scriptsize{4'}}
\put(4.04,.32){\scriptsize{5}}
\put(.45,.32){\scriptsize{-5}}
\put(4.45,.32){\scriptsize{6}}
\put(.04,.32){\scriptsize{-6}}

 \put(2.4,.4){\circle*{0.1}}
 \put(.4,.4){\circle*{0.1}}
\put(1.6,.4){\circle*{0.1}}
 \put(2,.4){\circle*{0.1}}
 \put(1.2,.4){\circle*{0.1}}
 \put(2.8,.4){\circle*{0.1}}
 \put(3.13,.4){\circle*{0.1}}\put(3.27,.4){\circle*{0.1}}
 \put(3.53,.4){\circle*{0.1}}\put(3.67,.4){\circle*{0.1}}
\put(4,.4){\circle*{0.1}}
\put(.8,.4){\circle*{0.1}}
 \put(4.4,.4){\circle*{0.1}}
 \put(0,.4){\circle*{0.1}}
\put(-1.25,.47){\vector(1,0){.8}}
\put(-1.,.56){$\psi$}
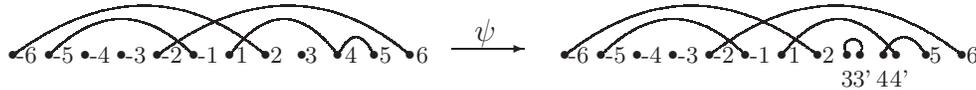
\put(-6.,-.5){\footnotesize{Fig. 4. The left digram is the upper permutation diagram of $\sigma$ = (4, -6, 3, 5, 1, -2) }}
\put(-5.,-.8){\footnotesize{and the right is its image by $\psi$. }}
\end{picture}\\\\
\textit{ Second step.} We give an outline of the involution $\varphi$ of \cite{KZ}. Let $\pi$ be a partition of type A and G be its partition diagram defined as the upper permutation diagram.
 For each two vertices $k$ and $j$ of G, we adopt that $j$ is a vacant vertex for the $k$th position if $j<k$ and its corresponding closer vertex $l$ satisfies $l$ $>$ $k$.
 Then for each arc (i, j) of G, we denote by $\delta(i,j)$ (resp. $\gamma(i,j)$) the number of vacant vertex $k$ such that $k < i$ (resp. $k>i$) for the position $j$.
The algorithm describing the involution $\varphi$ is to construct a partition diagram G$'$ from G, vertex by vertex and from left to right in the following paragraph.

For each vertex $k$ of G from 1 to the rank of $\pi$, if $k$ is a fixed (resp. opener) vertex then we conserve its form; fixed (resp. opener) vertex, at the position
 $k$ in G$'$ and  if $k$ is a closer (resp. transient) vertex, we also conserve its form; closer (resp. transient) vertex, at the same position in G$'$,
but we exchange the arc ($s$, $k$) where $s$ is the corresponding opener of $k$ in G
into an arc ($t$, $k$) in G$'$ with $t$ is the $\gamma(s,k)$th vacant vertex, from left to right, for the position $k$. So, $\varphi$ is a proper crossings and nestings interchanging map.\\
\begin{picture}(4,1.5)(-1,.2)

\qbezier(.4,.4)(1.2,1.2)(2,.4)
\qbezier(2.4,.4)(3.035,1.2)(3.67,.4)
\qbezier(3.51,.4)(3.755,.8)(4,.4)
\qbezier(0,.4)(1.4,1.5)(2.8,.4)
\qbezier(1.6,.4)(3.,1.5)(4.4,.4)

\qbezier(3.13,.4)(3.04,.55)(3.2,.58)
\qbezier(3.27,.4)(3.33,.55)(3.2,.58)

\put(1.25,.32){\scriptsize{-3}}
 \put(.85,.32){\scriptsize{-4}}
\put(1.65,.32){\scriptsize{-2}}
 \put(2.05,.32){\scriptsize{-1}}
 \put(2.475,.32){\scriptsize{1}}
 \put(2.85,.32){\scriptsize{2}}
 \put(3.06,.05){\scriptsize{3}}\put(3.19,.05){\scriptsize{3'}}
 \put(3.45,.05){\scriptsize{4}}\put(3.6,.05){\scriptsize{4'}}
\put(4.04,.32){\scriptsize{5}}
\put(.45,.32){\scriptsize{-5}}
\put(4.45,.32){\scriptsize{6}}
\put(.04,.32){\scriptsize{-6}}

\put(2.4,.4){\circle*{0.1}}
\put(.4,.4){\circle*{0.1}}
\put(1.6,.4){\circle*{0.1}}
\put(2,.4){\circle*{0.1}}
\put(1.2,.4){\circle*{0.1}}
\put(2.8,.4){\circle*{0.1}}
\put(3.13,.4){\circle*{0.1}}\put(3.27,.4){\circle*{0.1}}
\put(3.53,.4){\circle*{0.1}}\put(3.67,.4){\circle*{0.1}}
\put(4,.4){\circle*{0.1}}
\put(.8,.4){\circle*{0.1}}
\put(4.4,.4){\circle*{0.1}}
\put(0,.4){\circle*{0.1}}
\put(4.9,.54){\vector(1,0){.8}}
\put(5.25,.7){$\varphi$}
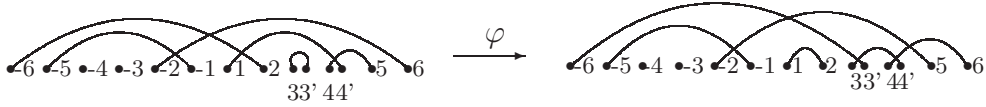
\put(2.2,-.7){\footnotesize{Fig. 5. The left diagram is the diagram of $\psi(\sigma)$ }}
\put(3.2,-1.){\footnotesize{and the right is its image by $\varphi$.}}
\end{picture}

\begin{picture}(4,.7)(-7.2,-.7)
\qbezier(.4,.4)(1.2,1.3)(2,.4)
\qbezier(1.6,.4)(2.8,1.6)(4,.4)
\qbezier(3.53,.4)(4,1)(4.4,.4)

\qbezier(2.4,.4)(2.6,.8)(2.8,.4)

\qbezier(3.13,.4)(3.4,.8)(3.67,.4)
\qbezier(0,.4)(1.6,1.8)(3.27,.4)

\put(1.25,.32){\scriptsize{-3}}
 \put(.8,.32){\scriptsize{-4}}
\put(1.64,.32){\scriptsize{-2}}
 \put(2.05,.32){\scriptsize{-1}}
 \put(2.45,.32){\scriptsize{1}}
 \put(2.82,.32){\scriptsize{2}}
 \put(3.1,.15){\scriptsize{3}}\put(3.24,.15){\scriptsize{3'}}
 \put(3.5,.15){\scriptsize{4}}\put(3.61,.15){\scriptsize{4'}}
\put(4.04,.32){\scriptsize{5}}
\put(.45,.32){\scriptsize{-5}}
\put(4.45,.32){\scriptsize{6}}
\put(.04,.32){\scriptsize{-6}}

 \put(2.4,.4){\circle*{0.1}}
 \put(.4,.4){\circle*{0.1}}
\put(1.6,.4){\circle*{0.1}}
 \put(2,.4){\circle*{0.1}}
 \put(1.2,.4){\circle*{0.1}}
 \put(2.8,.4){\circle*{0.1}}
 \put(3.12,.4){\circle*{0.1}}\put(3.26,.4){\circle*{0.1}}
 \put(3.52,.4){\circle*{0.1}}\put(3.66,.4){\circle*{0.1}}
\put(4,.4){\circle*{0.1}}
\put(.8,.4){\circle*{0.1}}
 \put(4.4,.4){\circle*{0.1}}
 \put(0,.4){\circle*{0.1}}
\end{picture}\\\\

\hspace{0.2cm}It remains to prove that the number of the minus signs is unchanged. Since, for each $\sigma \in B_n$,
the number of minus signs is equal to the number of arcs, in the upper permutation diagram $Upp(\sigma)$,
joining a vertex with negative index and a vertex with  positive index. Thus, according to the construction of the openers and closers vertices
of $\varphi(Upp(\sigma))$,
the last number is invariant.\\

Finally, for instance, by the map $\psi^{-1}\circ\varphi\circ\psi$ and Lemma 2. 1, we illustrate, in the following figure, the corresponding
permutation $\sigma'$ = (2, -5, 4, 6, 1, -3) of the permutation $\sigma$ = (4, $-$6, 3, 5, 1, $-$2) in Example 1.

\begin{picture}(1,3.2)(-4,-1.5)
\qbezier(.4,.4)(1.2,1.3)(2,.4)
\qbezier(1.6,.4)(2.8,1.6)(4,.4)\qbezier(1.2,.4)(2.6,-1.4)(4.4,.4)
\qbezier(3.6,.4)(4,1)(4.4,.4)
\qbezier(.4,.4)(1.6,-.9)(2.8,.4)
\qbezier(1.6,.4)(1.8,.)(2,.4)
\qbezier(2.4,.4)(2.6,.8)(2.8,.4)

\qbezier(2.4,.4)(3.2,-.6)(4,.4)
\qbezier(3.2,.4)(3.4,.8)(3.6,.4)
\qbezier(.8,.4)(1.,-.)(1.2,.4)

\qbezier(0,.4)(1.6,1.8)(3.2,.4)
\qbezier(0,.4)(0.4,-.15)(0.8,.4)

\put(1.25,.32){\scriptsize{-3}}
 \put(.8,.32){\scriptsize{-4}}
\put(1.64,.32){\scriptsize{-2}}
 \put(2.05,.32){\scriptsize{-1}}
 \put(2.45,.32){\scriptsize{1}}
 \put(2.82,.32){\scriptsize{2}}
 \put(3.27,.32){\scriptsize{3}}
 \put(3.65,.32){\scriptsize{4}}
\put(4.04,.32){\scriptsize{5}}
\put(.45,.32){\scriptsize{-5}}
\put(4.45,.32){\scriptsize{6}}
\put(.04,.32){\scriptsize{-6}}

 \put(2.4,.4){\circle*{0.1}}
 \put(.4,.4){\circle*{0.1}}
\put(1.6,.4){\circle*{0.1}}
 \put(2,.4){\circle*{0.1}}
 \put(1.2,.4){\circle*{0.1}}
 \put(2.8,.4){\circle*{0.1}}
 \put(3.2,.4){\circle*{0.1}}
 \put(3.6,.4){\circle*{0.1}}
\put(4,.4){\circle*{0.1}}
\put(.8,.4){\circle*{0.1}}
 \put(4.4,.4){\circle*{0.1}}
 \put(0,.4){\circle*{0.1}}
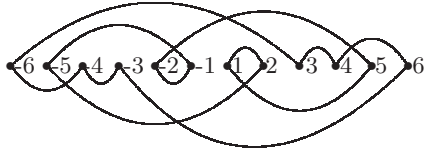
 \put(-1.2,-1){\footnotesize{Fig. 6. The permutation diagram of $\sigma$' = (2, -5, 4, 6, 1, -3).\hspace{2cm}$\Box$}}
\end{picture}
 \section{Proof of Theorem 2.7}
Our proof is based on an extension of de Mier's bijection in \cite[Section 4]{M} to $B_n$.\\
First, we extend the basic tool, in \cite{M}, that is the construction of a bijection between link partitions of type A and fillings of Young diagram into
a bijection, denoted by $\xi$, between upper permutation diagrams and fillings of Young diagrams on $\overline{B}_n$. For each $\sigma$ in $B_n$, let $i_1$, ..., $i_c$ be the closers
vertices of $Upp(\sigma)$ and $j_1$, ..., $j_o$ the openers ones. Let $p(i)$, for each closer vertex $i$, be the number of vertices $j$ with $j<i$ that are openers such
that for each transient we associate a closer before an opener. We consider a Young diagram $T$ of shape ($p(i_c)$, ..., $p(i_1)$), and if there is an arc going
from the opener $j_s$ to the closer $i_r$, we fill the cell in column $s$ and row $c - r + 1$ with~1. For instance, in the following figure, we
illustrate the upper permutation diagram of $\sigma$ of Example 2 and its corresponding filling of Young  diagram.

\begin{picture}(4,.)(-1.,-1.5)

\qbezier(0,.4)(1.,1.6)(2.4,.4)
\qbezier(.8,.4)(1.2,1.)(1.6,.4)

\put(1.25,.32){\scriptsize{-3}}
 \put(.8,.32){\scriptsize{-4}}
\put(1.64,.32){\scriptsize{-2}}
 \put(2.05,.32){\scriptsize{-1}}
 \put(2.45,.32){\scriptsize{1}}
 \put(2.82,.32){\scriptsize{2}}
 \put(3.27,.32){\scriptsize{3}}
 \put(3.65,.32){\scriptsize{4}}
\put(4.04,.32){\scriptsize{5}}
\put(.45,.32){\scriptsize{-5}}
\put(4.45,.32){\scriptsize{6}}
\put(.04,.32){\scriptsize{-6}}

 \put(2.4,.4){\circle*{0.07}}
 \put(.4,.4){\circle*{0.07}}
\put(1.6,.4){\circle*{0.07}}
 \put(2,.4){\circle*{0.07}}
 \put(1.2,.4){\circle*{0.07}}
 \put(2.8,.4){\circle*{0.07}}
 \put(3.2,.4){\circle*{0.07}}
 \put(3.6,.4){\circle*{0.07}}
\put(4,.4){\circle*{0.07}}
\put(.8,.4){\circle*{0.07}}
 \put(4.4,.4){\circle*{0.07}}
 \put(0,.4){\circle*{0.07}}
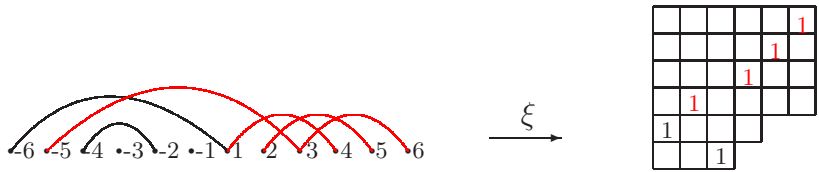
\put(.7,-.4){\footnotesize{Fig. 7. The upper permutation diagram of $\sigma$ = (4, 5, 6, 2, -3, -1) and }}
\put(1.5,-.8){\footnotesize{ its corresponding filling of Young  diagram.}}

\put(5.3,.54){\vector(1,0){.8}}
\put(5.65,.7){$\xi$}
\color{red}
\qbezier(.4,.4)(1.9,1.8)(3.2,.4)

\qbezier(2.4,.4)(3,1.2)(3.6,.4)
\qbezier(2.8,.4)(3.4,1.2)(4,.4)
\qbezier(3.2,.4)(3.8,1.2)(4.4,.4)

\end{picture}
\begin{picture}(4,4)(-8,-1.5)

\put(-4,2){\line(1,0){1.8}}
\put(-4,1.7){\line(1,0){1.8}}
\put(-4,1.4){\line(1,0){1.8}}
\put(-4,1.1){\line(1,0){1.8}}
\put(-4,.8){\line(1,0){1.8}}
\put(-4,.5){\line(1,0){1.2}}
\put(-4,.2){\line(1,0){.9}}

\put(-4,.2){\line(0,1){1.8}}
\put(-3.7,.2){\line(0,1){1.8}}
\put(-3.4,.2){\line(0,1){1.8}}
\put(-3.1,.2){\line(0,1){1.8}}
\put(-2.8,.5){\line(0,1){1.5}}
\put(-2.5,.8){\line(0,1){1.2}}
\put(-2.2,.8){\line(0,1){1.2}}
\put(-4,.54){\scriptsize{ 1}}
\put(-3.4,.25){\scriptsize{ 1}}
\color{red}
\put(-2.5,1.71){\scriptsize{ 1}}
\put(-2.8,1.42){\scriptsize{ 1}}
\put(-3.1,1.13){\scriptsize{ 1}}

\put(-3.7,.84){\scriptsize{ 1}}

\end{picture}\\
Thus, we see that the $k-$nesting (one can also say $k+1$-nonnesting) (resp., $k-$crossing) in the upper permutation diagram corresponds to a matrix identity $I_k$ (resp., antiidentity $J_k$ called also the antidiagonal of $I_k$) in a largest rectangle in the corresponding Young diagram.

\hspace{0.2cm}We can now state an important result due to de Mier \cite{M}.

\begin{lem} \cite[Theorem 3.5]{M} For all diagrams T with prescribed row and colunm sums, the number of fillings T that avoid $I_k$ equals the number of fillings of T that avoid $J_k$. \end{lem}
\hspace{0.2cm}Moreover, de Mier proves that there is a map $\Psi$ preserving the left-right degree
sequence of a link partition of type A and it interchanges the $k$-noncrossings and $k$-nonnestings for the proper crossings and nestings and conversely. We extend this involution on $B_n$. In fact, it divides into two maps. The first is $\varphi$.
For each $\sigma \in B_n$, we see the largest $k$ such that the filling of Young diagram corresponding to $Upp(\sigma)$ contains a largest rectangle which contains a matrix
antiidentity $J_k$. If there are many matrices antiidentities of rank $k$, we choose the one more to the right and the topmost. So, the 1's of $J_k$,
 from left and bottom to right and top, correspond to $(l_1,c_1)$, $(l_2, c_2)$, $\dots$, $(l_k, c_k)$ cells in the diagram, $i.e.$, $(l_i, c_i)$ is the
 intersection cell of the $l_i$th line and $c_i$th colunm, for each $1$  $\leqslant$ $i$ $\leqslant k$. Thus $\varphi$ changes the places of the 1's of $J_k$ in
 the diagram to new places define by:  $(l_2,c_1)$, $(l_3, c_2)$, $\dots$, $(l_k, c_{k-1})$ and $(l_1,c_k)$ and we obtain in the first time the following matrix
$$
\varphi(J_k) =  \begin{pmatrix}
 J_{k-1} & 0  \\
 0 & 1
\end{pmatrix}
$$ where $J_{k-1}$ is the matrix antiidentity of rank $k-1$.

So on, we apply $\varphi$ to  $J_{k-1}$, $J_{k-2}$, $\dots$, until we get $I_k$.\\\\
 \hspace{0.2cm}The second is $\phi$ (the inverse of $\varphi$). Let $k$ be the largest integer such that the Young diagram contains a matrix $I_k$
(if there are many matrices identities of rank $k$, we choose the one more to the right and the topmost). The 1's of $I_k$ have $(a_1,b_1)$, $(a_2, b_2)$, $\dots$, $(a_k, b_k)$ as the cells in the diagram
 from top to bottom and left to right with $(a_i, b_i)$ is the intersection cell between the $a_i$th line and the $b_i$ colunm. So, the image of $I_k$ by $\phi$ is
$(a_2,b_1)$, $(a_1, b_2)$, $(a_3, b_3)$, $\dots$, $(a_k, b_k)$, $i.e.$,
 $$\phi(I_k) =  \begin{pmatrix}
 0 & 1 & 0  \\
 1 & 0 & 0 \\
 0 & 0 & I_{k-2}
\end{pmatrix} $$
 where $I_{k-2}$ is the matrix identity of rank $k-2$. The image of $\phi(I_k)$ by $\phi$ is $(a_3,b_1)$, $(a_2, b_2)$, $(a_1, b_3)$, $(a_4, b_4)$, $\dots$, $(a_k, b_k)$.
 So on, we apply $\phi$ to  $\phi^i(I_k)$, for $i$ from 0 to $k-1$, until we get $J_k$.

 Then, by the two processes, $\Psi$ interchanges $I_k$ and $J_k$.

\hspace{0.2cm}It remains to prove that the number of minus signs is unchanged by the above transformations $\xi$ and $\Psi$. We know that this number in such a permutation
 in $B_n$ is equal to the number of the arcs ($\alpha_i$, $\beta_i$) in its upper permutation diagram such that $\alpha_i$ (resp., $\beta_i$) is a negative (resp., positive) index.
 Let $\sigma$ in $B_n$. Suppose that $neg(\sigma)$ = m. Then there exist only $m$ arcs ($\alpha_1$, $\beta_1$), \ldots, ($\alpha_m$, $\beta_m$) in
 $Upp(\sigma)$ such that, for each $1 \leq i \leq m$, $\alpha_i$ and $\beta_i$ are, respectively, a negative and positive indices.
 Let ($i$, $j$) be an arc in $Upp(\sigma)$. We say that $i$ (resp., $j$) is an $outpoint$ (resp., $endponit$) and we denote by $NO$ (resp., $NE$) the number of outpoints (resp., endpoints) with negative
index. We have that $m$ = ($NO$ - $NE$). On the other hand, the degree sequence of $Upp(\sigma)$ is kept by the involution $\vartheta$ where
$\vartheta$ = $\xi^{-1}\circ\Psi\circ\xi$. Then there exist only ($NO$ - $NE$) arcs in $\vartheta(Upp(\sigma))$ such that each arc has an outpoint with negative
 index and an endpoint with positive index. This establishes the desired equality.

For instance, the following diagram is the permutation diagram of $\vartheta(Upp(\sigma))$ of the permutation $\sigma$ in Example 2.
 It is easy to see that $cro^*(\vartheta(Upp(\sigma)))$ = 2 and $nes^*(\vartheta(Upp(\sigma)))$ = 4.\\

\begin{picture}(9.,2.3)(-7,-.8)
\qbezier(.4,.4)(1,1.2)(1.6,.4)
\qbezier(0,.4)(1.,1.6)(2.4,.4)

\put(1.25,.32){\small{-3}}
 \put(.8,.32){\small{-4}}
\put(1.64,.32){\small{-2}}
 \put(2.05,.32){\small{-1}}
 \put(2.45,.32){\small{1}}
 \put(2.82,.32){\small{2}}
 \put(3.27,.32){\small{3}}
 \put(3.65,.32){\small{4}}
\put(4.04,.32){\small{5}}
\put(.45,.32){\small{-5}}
\put(4.45,.32){\small{6}}
\put(.04,.32){\small{-6}}

 \put(2.4,.4){\circle*{0.07}}
 \put(.4,.4){\circle*{0.07}}
\put(1.6,.4){\circle*{0.07}}
 \put(2,.4){\circle*{0.07}}
 \put(1.2,.4){\circle*{0.07}}
 \put(2.8,.4){\circle*{0.07}}
 \put(3.2,.4){\circle*{0.07}}
 \put(3.6,.4){\circle*{0.07}}
\put(4,.4){\circle*{0.07}}
\put(.8,.4){\circle*{0.07}}
 \put(4.4,.4){\circle*{0.07}}
 \put(0,.4){\circle*{0.07}}
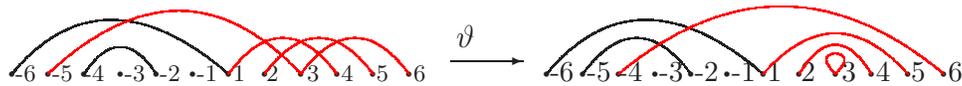
 \put(-4.2,-.5){\footnotesize{Fig. 8. The upper permutation diagram $\vartheta(Upp(\sigma))$}.\hspace{3.5cm}$\Box$}
\put(-1.06,.54){\vector(1,0){.8}}
\put(-1.,.7){$\vartheta$}
\color{red}
\qbezier(2.4,.4)(3.2,1.3)(4,.4)
\qbezier(2.8,.4)(3.2,1.05)(3.6,.4)
\qbezier(3.2,.4)(3,.6)(3.2,.63)\qbezier(3.2,.4)(3.4,.6)(3.2,.63)
\qbezier(.8,.4)(2.6,1.9)(4.4,.4)

\end{picture}
\begin{picture}(4,1.5)(8.,-.8)
\qbezier(0,.4)(1.,1.6)(2.4,.4)
\qbezier(.8,.4)(1.2,1.)(1.6,.4)

\put(1.25,.32){\scriptsize{-3}}
 \put(.8,.32){\scriptsize{-4}}
\put(1.64,.32){\scriptsize{-2}}
 \put(2.05,.32){\scriptsize{-1}}
 \put(2.45,.32){\scriptsize{1}}
 \put(2.82,.32){\scriptsize{2}}
 \put(3.27,.32){\scriptsize{3}}
 \put(3.65,.32){\scriptsize{4}}
\put(4.04,.32){\scriptsize{5}}
\put(.45,.32){\scriptsize{-5}}
\put(4.45,.32){\scriptsize{6}}
\put(.04,.32){\scriptsize{-6}}

 \put(2.4,.4){\circle*{0.07}}
 \put(.4,.4){\circle*{0.07}}
\put(1.6,.4){\circle*{0.07}}
 \put(2,.4){\circle*{0.07}}
 \put(1.2,.4){\circle*{0.07}}
 \put(2.8,.4){\circle*{0.07}}
 \put(3.2,.4){\circle*{0.07}}
 \put(3.6,.4){\circle*{0.07}}
\put(4,.4){\circle*{0.07}}
\put(.8,.4){\circle*{0.07}}
 \put(4.4,.4){\circle*{0.07}}
 \put(0,.4){\circle*{0.07}}

\color{red}
\qbezier(.4,.4)(1.9,1.8)(3.2,.4)

\qbezier(2.4,.4)(3,1.2)(3.6,.4)
\qbezier(2.8,.4)(3.4,1.2)(4,.4)
\qbezier(3.2,.4)(3.8,1.2)(4.4,.4)

\end{picture}\\

\hspace{0.2cm}We now come to conclude this paper with a straightforward enumerative result on the maximum crossing (resp. maximum nesting) chains.
 A $maximum$ $crossing$, in such a permutation in $B_n$, is the $n$-crossing. Then we can compute the number of permutations with $n$-crossing, denoted by $C_{B_n}$,
in the following corollary.\\

\begin{cor}
\textit{Let n be a positive integer. Then the number of permutations in $B_n$ with n-crossings $C_{B_n}$ satisfies }
\begin{align}\label{1}
C_{B_n} = \left\{ \begin{array}{lll}
 2\;\,  & & if\;\  n \leq 2,\\
 1 & & otherwise.
 \end{array}\right.
\end{align}
\end{cor}
In fact, for $n \leq 2$, it is easy to see that. For $n \geq 3$, the only element $\sigma$ in $B_n$ that has n-crossing is the permutation $\sigma$ defined by
 $\sigma$($i$) = $-$($n + 1$) + $i$ for $1 \leq i \leq n$, $i.e.$, the only permutation $\sigma$ that satisfies $\sigma(1)<\sigma(2)< \ldots<\sigma(n)<1$.

\section*{Acknowledgements}
\hspace{0.2cm}The author would like to thank Jiang Zeng for his advice and careful reading of this paper.



\begin{thebibliography}{99}
\bibitem{BMP} S. Burrill, M. Mishna and J. Post, On  $k$-crossings and $k$-nestings of permutations, DMTCS proc. AN, (2010), 461-468.
 \bibitem{CDDS} W. Y. C. Chen, E. Y. P. Deng, R. R. X. Du, R. P. Stanley and C. H. Yan, Crossings and nestings of matchings and partitions, Trans. Amer. Math. Soc. 359(4)(2007), 1555-1575.
 \bibitem{Co}  S. Corteel, Crossings and alignnments of permutations, Adv. in Appl. Math., 38(2)(2007), 149-163.
   \bibitem{CMW}  S. Corteel, M. Josuat-Verg\`es,  L. Williams, The matrix ansatz, orthogonal polynomials and permutations, ariXiv:math.CO/1005.2696v1, to appear in Adv. Appl. Math.
  \bibitem{KZ}  A. Kasraoui, J. Zeng, Distrbution of crosings, nestings and alignments of two edges in matchings and partitions, Electron. J. Combin. 13(1) (2006), R33.
 \bibitem{Kr} C. Krattenthaler, Growth diagrams, and increasing and decrasing chains in fillings of Ferrers sharps, Adv. in Appl. Math. 37(3)(2006), 404-431.
 \bibitem{M} A. De Mier, k-noncrossing and k-nonnesting graphs and fillings of Ferrers diagrams, Combinatorica  27(6) (2007), 699-720.
\bibitem{RS} M. Rubey and C. Stump, Crossings and nestings in set partitions of classical types, Electron. J.  Combin. 17(1) (2010), R120.
\end{thebibliography}
\end{document}